# Divisor Problem and an Analogue of Euler's Summation Formula

VIVEK V. RANE

Department of Mathematics
The Institute of Science,
15, Madam Cama Road,
Mumbai-400 032
India .
v_v_rane@yahoo.co.in

Abstract : By elementary way , we obtain with ease , a highly simple expression for $\Delta(x)$, the remainder term of Dirichlet's divisor problem and use it to derive an analogue of Euler's summation formula for the summation $\sum_{a < n \leq b} d(n) f(n)$ , where $d(n)$ is the divisor function and $f(x)$ is a smooth function on [a,b] . For the Riemann zeta function $\zeta(s)$ with $0 < s \neq 1$ and $x > 0$ , if $\zeta(s) = \sum_{n \leq x} n^{-s} + E_1(s,x)$ and $\zeta^2(s) = \sum_{n \leq x} d(n) n^{-s} + E_2(s,x)$ , we state a concise relation among $E_1^2(s,x), E_2(s,x)$ and $\Delta(x)$ , while correcting an error in an earlier paper ( of the author) , dealing with divisor problem .

# Divisor Problem and an Analogue of Euler's Summation Formula

## VIVEK V. RANE

**Introduction :** For a natural number n, let d(n) denote the number of divisors of n. and let $D(x) = \sum_{n \leq x} d(n)$. Let $D(x) = x \log x + (2\gamma - 1)x + \frac{1}{4} + \Delta(x)$, where $\gamma$ is Euler's constant. Our object is to obtain a simple expression for $\Delta(x)$ by elementary method, especially when $x$ is a positive integer and use it to obtain an analogue of Euler-Poisson summation formula for the summation $\sum_{a < n \leq b} d(n) f(n)$, where $f$ is a continuously differentiable function on the interval [a,b] with $b > a > 0$. We shall write $f'$ to denote its derivative. In what follows, for a real variable u, we shall write [u] to denote its integral part. We write $\psi(u) = u - [u] - \frac{1}{2}$. For any integer $r \geq 1$, we write

$$\psi_r(u) = -\sum_{|n| \geq 1} \frac{e^{2\pi i n u}}{(2\pi i n)^r}$$ so that $\psi_r(u) = \frac{B_r(u-[u])}{r!}$, where $B_r(u)$ is Bernoulli polynomial of degree r. Note $\psi_1(u) = \psi(u)$ and $\frac{d}{du}\psi_r(u) = \psi_{r-1}(u)$ for $r \geq 1$. In what follows, O-constants are absolute.

Next, we state our Theorems.

**Theorem 1 :** Let $D(x) = \sum_{n \leq x} d(n) = x \log x + (2\gamma - 1)x + \frac{1}{4} + \Delta(x)$, where $\gamma$ is Euler's constant and let $\psi(u) = u - [u] - \frac{1}{2}$. Then we have $\Delta(x) = A(x) + B(x)$, where

$A(x) = -2 \sum_{m \leq \sqrt{x}} \psi(\frac{x}{m})$ and $B(x) = 4x \cdot \int_{\sqrt{x}}^{\infty} \frac{\psi_2(u)}{u^3} du - \psi^2(\sqrt{x}) - 2\psi_2(\sqrt{x}) = 0(1)$.

Here $\psi_2(x) = \frac{B_2(x-[x])}{2}$, where $B_2(x)$ is Bernoulli polynomial of degree 2.

In addition, if $x$ is an integer, $A(x) = \frac{1}{\pi i} \sum_{m \leq \sqrt{x}} \sum_{r=1}^{m-1} e^{\frac{2\pi i r x}{m}} \left( \frac{1}{r} - 2r \sum_{v \geq 1} \frac{1}{m^2 v^2 - r^2} \right) + 0(d(x))$.

**Theorem 2 :** Let $f$ be a continuously differentiable function on an interval [a,b], where, $b > a > 0$. Then we have

$$\sum_{a<n\leq b} d(n)f(n) = \int_a^b (\log u + 2\gamma)f(u)du - 2\int_a^b \frac{f(u)\psi_2(\sqrt{u})}{u}du + \int_a^b \frac{f(u)\psi(\sqrt{u})}{\sqrt{u}}du$$

$$+ 4\int_a^b f(u)\left(\int_{\sqrt{u}}^\infty \frac{\psi_2(t)}{t^3}dt\right)du - f(b)\psi^2(\sqrt{b}) + f(a)\psi^2(\sqrt{a})$$

$$+ \int_a^b \psi^2(\sqrt{u})f'(u)du + 2f(a)\sum_{m\leq\sqrt{a}}\psi(\tfrac{a}{m}) - 2f(b)\sum_{m\leq\sqrt{b}}\psi(\tfrac{b}{m}) - \tfrac{1}{\pi i}\sum_{m\leq\sqrt{b}}\sum_{|v|\geq 1}\tfrac{1}{v}\cdot\int_{\max(m^2,a)}^b e^{\frac{2\pi i v u}{m}}\cdot f'(u)du$$

**Remark** : Note that $\int_{\sqrt{u}}^\infty \frac{\psi_2(t)}{t^3}dt = O(\tfrac{1}{u})$.

It is to be noted that the relation $\Delta(u) = -2\sum_{m\leq\sqrt{u}}\psi(\tfrac{u}{m}) + O(1)$ is well-known. However, we have given highly simple explicit elementary expressions for $A(u)$ and $B(u)$, where $\Delta(u) = A(u) + B(u)$. In author [2], we have given analogues of Euler and Poisson summation formulae for the summation $\sum_{a<n\leq b} d(n)f(n)$. Our Theorem 2 gives a somewhat different analogue of Euler-Poisson summation formula for $\sum_{a<n\leq b} d(n)f(n)$, based on Divisor Problem. However, in principle, these two analogues for $\sum_{a<n\leq b} d(n)f(n)$ are not very much different.

For $0 < s \neq 1$ and for $x > 0$, let $\zeta(s) = \sum_{n\leq x} n^{-s} + E_1(s,x)$ and let $\zeta^2(s) = \sum_{n\leq x} d(n)n^{-s} + E_2(s,x)$, where $\zeta(s)$ is Riemann zeta function. These are approximate functional equations for $\zeta(s)$ and $\zeta^2(s)$ respectively for real s.

In author [3], we have shown that $E_2(s,x) = 2\sum_{n\leq\sqrt{x}} n^{-s} E_1(s,\tfrac{x}{n}) + E_1^2(s,\sqrt{x})$ and $E_1(s,x) = \tfrac{x^{1-s}}{s-1} + x^{-s}\psi(x) + O(sx^{-s-1})$ for $s > 0$ and $s \neq 1$. In the light of the fact that

$$\Delta(x) = -2\sum_{n\leq\sqrt{x}}\psi(\tfrac{x}{n}) + O(1),$$

we have the following Theorem 3.

**Theorem 3** : For $0 < s \neq 1$ and for $x > 0$,

let $\zeta(s) = \sum_{n \leq x} n^{-s} + E_1(s,x)$ and let $\zeta^2(s) = \sum_{n \leq x} d(n) n^{-s} + E_2(s,x)$.

Then we have, $E_2(s,x) = \frac{2x^{1-s}}{s-1} \sum_{n \leq \sqrt{x}} \frac{1}{n} - x^{-s} \Delta(x) + E_1^2(s, \sqrt{x}) + 0(sx^{-s})$.

Note that in Theorem 2 of author [3], we have incorrectly used the fact that

$\sum_{n \leq x} \frac{1}{n} = \log x + \gamma + \frac{\gamma_1}{x} + \frac{\gamma_2}{x^2} + 0\left(\frac{1}{x^3}\right)$, where $\gamma, \gamma_1, \gamma_2$ are constants. This is true for integral $x$ only and <u>not</u> for general $x$. Instead, we write for general $x > 0$,

$\sum_{n \leq x} \frac{1}{n} = \log x + \gamma + \frac{\gamma_1(x)}{x} + \frac{\gamma_2(x)}{x^2} + 0\left(\frac{1}{x^3}\right)$, where $\gamma_1(x), \gamma_2(x)$ are functions of $x$. Note that $\gamma_1(x), \gamma_2(x)$ are constants, when $x$ is integral.

This gives $\sum_{n \leq \sqrt{x}} \frac{1}{n} = \log \sqrt{x} + \gamma + \frac{\gamma_1(\sqrt{x})}{\sqrt{x}} + \frac{\gamma_2(\sqrt{x})}{x} + 0\left(\frac{1}{x^{3/2}}\right)$. Then the statement of Theorem 2 of author [3] is true with $\gamma_1$ replaced by $\gamma_1(\sqrt{x})$ and $\gamma_2$ replaced by $\gamma_2(\sqrt{x})$ after noting that $\gamma_1(\sqrt{x}) + \psi(\sqrt{x}) = 0$. Our Theorem 3 above is a concise form of Theorem 2 of author [3].

Next, we give our proof of the Theorems.

**Proof of Theorem 1**: We have $D(x) = 2 \sum_{m \leq \sqrt{x}} \sum_{n \leq \frac{x}{m}} 1 - \sum_{m, n \leq \sqrt{x}} 1$.

Next $2 \sum_{n \leq \frac{x}{m}} 1 = 2\frac{x}{m} - 1 - 2\psi\left(\frac{x}{m}\right)$.

This gives $2 \sum_{m \leq \sqrt{x}} \sum_{n \leq \frac{x}{m}} 1 = 2x \sum_{m \leq \sqrt{x}} \frac{1}{m} - [\sqrt{x}] - 2 \sum_{m \leq \sqrt{x}} \psi\left(\frac{x}{m}\right)$.

Next $2x \sum_{m \leq \sqrt{x}} \frac{1}{m} = x \log x - 2x \int_1^{\sqrt{x}} \frac{\psi(u)}{u^2} du + x - 2\sqrt{x} \cdot \psi(\sqrt{x})$,

after using Euler's summation formula.

Next $-\int_1^{\sqrt{x}} \frac{\psi(u)}{u^2} du = -\int_1^{\infty} \frac{\psi(u)}{u^2} du + \int_{\sqrt{x}}^{\infty} \frac{\psi(u)}{u^2} du = \gamma - \frac{1}{2} + \int_{\sqrt{x}}^{\infty} \frac{\psi(u)}{u^2} du$.

Now, $\int_{\sqrt{x}}^{\infty} \frac{\psi(u)}{u^2} du = \int_{\sqrt{x}}^{\infty} \frac{d\psi_2(u)}{u^2} du = \left[\frac{\psi_2(u)}{u^2}\right]_{u=\sqrt{x}}^{\infty} + 2 \int_{\sqrt{x}}^{\infty} \frac{\psi_2(u)}{u^3} du = -\frac{\psi_2(\sqrt{x})}{x} + 2 \int_{\sqrt{x}}^{\infty} \frac{\psi_2(u)}{u^3} du$.

Note that $\int_{\sqrt{x}}^{\infty} \frac{\psi_2(u)}{u^3} du = 0(\frac{1}{x})$. Thus $2 \sum_{m \leq \sqrt{x}} \sum_{n \leq \frac{x}{m}} 1$

$= x \log x + (2\gamma - 1)x - 2\psi_2(\sqrt{x}) + 4x \int_{\sqrt{x}}^{\infty} \frac{\psi_2(u)}{u^3} du + x - [\sqrt{x}] - 2\sqrt{x} \cdot \psi(\sqrt{x}) - 2 \sum_{m \leq \sqrt{x}} \psi\left(\frac{x}{m}\right)$.

Next $\sum_{m, n \leq \sqrt{x}} 1 = [\sqrt{x}]^2 = \left(\sqrt{x} - \psi(\sqrt{x}) - \frac{1}{2}\right)^2 = x + \psi^2(\sqrt{x}) + \frac{1}{4} - 2\sqrt{x}\psi(\sqrt{x}) + \psi(\sqrt{x}) - \sqrt{x}$

Hence

$2 \sum_{m \leq \sqrt{x}} \sum_{n \leq \frac{x}{m}} 1 - \sum_{m, n \leq \sqrt{x}} 1 = x \log x + (2\gamma - 1)x - 2\psi_2(\sqrt{x}) + 4x \int_{\sqrt{x}}^{\infty} \frac{\psi_2(u)}{u^3} du + \frac{1}{4} - 2 \sum_{m \leq \sqrt{x}} \psi\left(\frac{x}{m}\right) - \psi^2(\sqrt{x})$

$= x \log x + (2\gamma - 1)x + \frac{1}{4} - 2 \sum_{m \leq \sqrt{x}} \psi\left(\frac{x}{m}\right) + B(x)$,

where $B(x) = 4x \int_{\sqrt{x}}^{\infty} \frac{\psi_2(u)}{u^3} du - 2\psi_2(\sqrt{x}) - \psi^2(\sqrt{x}) = 0(1)$.

This gives $\Delta(x) = -2 \sum_{m \leq \sqrt{x}} \psi\left(\frac{x}{m}\right) + B(x)$, where $B(x) = 0(1)$.

Next, we shall obtain an expression for $A(x) = -2 \sum_{m \leq \sqrt{x}} \psi\left(\frac{x}{m}\right)$, when $x$ is an integer.

For integral $x$, if $m$ does not divide $x$, we have

$-2\psi\left(\frac{x}{m}\right) = \sum_{|v| \geq 1} \frac{e^{2\pi i v \frac{x}{m}}}{\pi i v} = \frac{1}{\pi i} \sum_{r=1}^{m} e^{2\pi i r \frac{x}{m}} \cdot \sum_{\substack{|v| \geq 1 \\ v \equiv r \pmod{m}}} \frac{1}{v}$

Note that $\sum_{\substack{|v| \geq 1 \\ v \equiv m \pmod{m}}} \frac{1}{v} = \sum_{\substack{|v| \geq 1 \\ v \equiv 0 \pmod{m}}} \frac{1}{v} = 0$.

For $r = 1, 2, \ldots, m-1$,

we have $\sum_{v \equiv r \pmod{m}} \frac{1}{v} = \frac{1}{r} + \sum_{v \geq 1} \left(\frac{1}{-mv+r} + \frac{1}{mv+r}\right) = \frac{1}{r} + \sum_{v \geq 1} \frac{2r}{r^2 - m^2 v^2} = \frac{1}{r} - 2r \sum_{v \geq 1} \frac{1}{m^2 v^2 - r^2}$

Thus if $m$ does not divide $x$, then

$-2\psi(\frac{x}{m}) = \frac{1}{\pi i} \sum_{r=1}^{m-1} e^{2\pi i r \frac{x}{m}} \left(\frac{1}{r} - 2r \sum_{v \geq 1} \frac{1}{m^2 v^2 - r^2}\right)$.

If $m$ divides $x$, then $-2\psi(\frac{x}{m}) = -2(\frac{x}{m} - [\frac{x}{m}] - \frac{1}{2}) = (-2)(-\frac{1}{2}) = 1$ and $(-2) \sum_{|v| \geq 1} \frac{e^{2\pi i v \frac{x}{m}}}{-2\pi i v} = 0$.

This gives $-2 \sum_{m \leq \sqrt{x}} \psi(\frac{x}{m}) = \frac{1}{\pi i} \sum_{m \leq \sqrt{x}} \sum_{r=1}^{m-1} e^{\frac{2\pi i r x}{m}} \left(\frac{1}{r} - \sum_{v \geq 1} \frac{1}{m^2 v^2 - r^2}\right) + 0(d(x))$.

This completes the proof of Theorem 1.

Before proving Theorem 2, we state a few results concerning Riemann-Stieltjes integration in the form of a lemma. For the theory of Riemann-Stieltjes Integration, we refer the reader to Apostol's book [1]. However, our Lemma below gives a few additional results also.

**Lemma**: Let $f$ be a continuous function on an interval [a,b].

I) Let $g$ be a step function on [a, b]. Then $\int_a^b f dg$ exists.

II) Let $g$ be a function on [a,b] with its derivative $g'$ Riemann integrable on [a,b]. Then $\int_a^b f(u) dg(u) = \int_a^b f(u) g'(u) du$.

III) Let $G(u) = \int_a^u g(t) dt$, where $g$ is continuous on [a,b] except for finitely many simple discontinuities on [a,b]. Then $\int_a^b f(u) dG(u) = \int_a^b f(u) g(u) du$

IV) Let $g$ be a continuously differentiable function on [a, b] and let $h$ be a step function on [a,b]. Then $\int_a^b f(u) d(g(u) \cdot h(u)) = \int_a^b f(u) g(u) dh(u) + \int_a^b f(u) \cdot g'(u) h(u) du$.

**Proof of Theorem 2**: Let $f$ be a continuously differentiable function on [a,b] and let $D(u) = \sum_{n \le u} d(n)$. Then $\sum_{a < n \le b} d(n) f(n) = \int_a^b f(u) dD(u)$

Noting that $D(u) = u \log u + \frac{1}{4} - 2u \int_1^{\sqrt{u}} \frac{\psi(t)}{t^2} dt - \psi^2(\sqrt{u}) + A(u)$, where $A(u) = -2 \sum_{n \le \sqrt{u}} \psi\left(\frac{u}{n}\right)$,

we have $\sum_{a < n \le b} d(n) f(n) = \int_a^b f(u) dD(u)$

$= \int_a^b f(u) d\left( u \log u + \frac{1}{4} - 2u \cdot \int_1^{\sqrt{u}} \frac{\psi(t)}{t^2} dt - \psi^2(\sqrt{u}) + A(u) \right)$.

We shall see that each of the integrals $\int_a^b f(u) d(u \log u + \frac{1}{4})$ $\left( = \int_a^b f(u) d(u \log u) \right)$,

$\int_a^b f(u) d\left( -2u \cdot \int_1^{\sqrt{u}} \frac{\psi(t)}{t^2} dt \right)$, $\int_a^b f(u) d\psi^2(\sqrt{u})$ exists so that $\int_a^b f(u) dA(u)$ exists.

Now $\int_a^b f(u)d(u\log u + \frac{1}{4}) = \int_a^b f(u)d(u\log u) = \int_a^b f(u)(1+\log u)du$ .

Next $\int_a^b f(u)d\left(-2u \cdot \int_1^{\sqrt{u}} \frac{\psi(t)}{t^2}dt\right) = -2\int_a^b f(u)d\left(u\int_1^{\sqrt{u}} \frac{\psi(t)}{t^2}dt\right)$

$= -2\left(\int_a^b f(u) \cdot \int_1^{\sqrt{u}} \frac{\psi(t)}{t^2}dt + \int_a^b f(u)u \cdot \frac{d}{du}\left(\int_1^{\sqrt{u}} \frac{\psi(t)}{t^2}dt\right)\right)$

$= \int_a^b f(u)\left(-2\int_1^{\sqrt{u}} \frac{\psi(t)}{t^2}dt\right) - 2\int_a^b f(u) \cdot u \cdot \frac{\psi(\sqrt{u})}{(\sqrt{u})^2} \cdot \left(-\frac{1}{2\sqrt{u}}\right)du$

$= \int_a^b f(u)\left((2\gamma-1) - 2\frac{\psi_2(\sqrt{u})}{u} + 4\int_{\sqrt{u}}^{\infty} \frac{\psi_2(t)}{t^3}dt\right)du + \int_a^b \frac{f(u)\psi(\sqrt{u})}{\sqrt{u}}du$

$= (2\gamma-1)\int_a^b f(u)du - 2\int_a^b \frac{f(u)\psi_2(\sqrt{u})}{u}du + 4\int_a^b f(u)\left(\int_{\sqrt{u}}^{\infty} \frac{\psi_2(t)}{t^3}dt\right)du + \int_a^b \frac{f(u)\psi(\sqrt{u})}{\sqrt{u}}du$

Thus $\int_a^b f(u)d(u\log u) + \int_a^b f(u)d\left(-2u\int_1^{\sqrt{u}} \frac{\psi(t)}{t^2}dt\right)$

$= \int_a^b f(u)(\log u + 2\gamma)du - 2\int_a^b \frac{f(u)\psi_2(\sqrt{u})}{u}du + 4\int_a^b f(u)\left(\int_{\sqrt{u}}^{\infty} \frac{\psi_2(t)}{t^3}dt\right)du + \int_a^b \frac{f(u)\psi(\sqrt{u})}{\sqrt{u}}du$ ,

where $\int_{\sqrt{u}}^{\infty} \frac{\psi_2(t)}{t^3}dt = 0(\frac{1}{u})$ .

Next consider $\int_a^b f(u)d\psi^2(\sqrt{u})$ . Note that $\psi(u) = u - [u] - \frac{1}{2}$ .

This gives $\psi^2(\sqrt{u}) = u + [\sqrt{u}]^2 + \frac{1}{4} - 2\sqrt{u} \cdot [\sqrt{u}] + [\sqrt{u}] - \sqrt{u}$ .

Note that $[\sqrt{u}]$ is a step function . Hence in view of Lemma I) and IV) , $\int_a^b f(u)d\psi^2(\sqrt{u})$

exists . Also note that $-\int_a^b f(u)d\psi^2(\sqrt{u}) = -\left\{f(u)\psi^2(\sqrt{u})\big|_{u=a}^b - \int_a^b \psi^2(\sqrt{u})f'(u)du\right\}$

$= f(a)\psi^2(\sqrt{a}) - f(b)\psi^2(\sqrt{b}) + \int_a^b \psi^2(\sqrt{u})f'(u)du$ .

This gives $\int_a^b f(u)d\left(u\log u - 2u \cdot \int_1^{\sqrt{u}} \frac{\psi(t)}{t^2}dt + \frac{1}{4} - \psi^2(\sqrt{u})\right)$

: 7 :

$$= \int_a^b f(u)(\log u + 2\gamma)du - 2\int_a^b \frac{f(u)\psi_2(\sqrt{u})}{u}du + 4\int_a^b f(u)\left(\int_{\sqrt{u}}^{\infty} \frac{\psi_2(t)}{t^3}dt\right)du + \int_a^b \frac{f(u)\psi(\sqrt{u})}{\sqrt{u}}du$$

$$+ \int_a^b \psi^2(\sqrt{u})f'(u)du + f(a)\psi^2(\sqrt{a}) - f(b)\psi^2(\sqrt{b}) \ .$$

Next $\int_a^b f(u)dA(u) = f(u)A(u)|_{u=a}^b - \int_a^b A(u)f'(u)du$

$$= f(b)A(b) - f(a) \cdot A(a) - \int_a^b A(u)f'(u)du$$

$$= f(b)A(b) - f(a)A(a) + 2\int_a^b \left(\sum_{m \le \sqrt{u}} \psi(\tfrac{u}{m})\right) f'(u)du$$

$$= f(b)A(b) - f(a)A(a) + 2\sum_{m \le \sqrt{b}} \int_{\max(a,m^2)}^b \psi(\tfrac{u}{m}) \cdot f'(u)du$$

$$= f(b)A(b) - f(a)A(a) + 2\sum_{m \le \sqrt{b}} \int_{\max(a,m^2)}^b \left(-\sum_{|v| \ge 1} \frac{e^{2\pi i v \frac{u}{m}}}{2\pi i v}\right) f'(u)du$$

$$= f(b)A(b) - f(a)A(a) - \tfrac{1}{\pi i} \sum_{m \le \sqrt{b}} \sum_{|v| \ge 1} \tfrac{1}{v} \int_{\max(a,m^2)}^b e^{2\pi i v \frac{u}{m}} f'(u)du \ .$$

This completes the proof of Theorem 2 .